\documentclass[preprint,1p]{elsarticle}
\usepackage{graphics}
\usepackage{epsfig}
\usepackage{mathptmx}
\usepackage{amsmath}
\usepackage{amssymb}
\usepackage{hyperref}
\usepackage{fullpage}
\usepackage{amsthm}
\usepackage{lineno}
\journal{}
\begin{document}
\begin{frontmatter}
\title{Kantorovich form of generalized Szasz-type operators with certain parameters using Charlier polynomials}
\author[]{Abdul Wafi}
\ead{awafi@jmi.ac.in}
\author[]{Nadeem Rao\corref{cor1}}
\ead{nadeemrao1990@gmail.com}
\address{Department of Mathematics, Jamia Millia Islamia, New Delhi-110 025, India}
\cortext[cor1]{Corresponding author}
\begin{abstract}

The aim of this article is to introduce the Kantorovich form of generalized Szasz-type operators involving Charlier polynomials with certain parameters. In this paper we discussed the rate of convergence, better error estimates and Korovkin-type theorem in polynomial weighted space. Further, we investigate the local approximation results with the help of Ditzian-Totik modulus of smoothness, second order modulus of continuity, Peetre's K-functional and Lipschitz class.
\end{abstract}
\begin{keyword}

Szasz operators, Charlier polynomials, Ditzian-Totik modulus of smoothness, Peetre's K-functional, Lipschitz class.
\vspace{0.5 cm}
\newline\textbf{2010 Mathematics Subject Classification 41A10, 41A25, 41A36, 41A36}
\end{keyword}
\end{frontmatter}

\section{Introduction}

Approximation theory plays an important role in mathematical analysis and other branches of mathematics. The results of theory of approximation is generally related to positive linear operators, and deals with rate of convergence and order of approximation. Weierstrass was the first who gave an important theorem, namely, Weierstrass approximation theorem in this regard. The aim of this theorem is to minimize the maximum value of $|f(x)-P_n(x)|$ for the continuous functions $f(x)$ on $[a,b]$, where $P_n(x)$ is the polynomial of degree $n$. The proof of this theorem was considered very difficult until Bernstein gave an elegant and simple proof of it. Beernstein\cite{Bernstein} defined the positive linear operators using binomial distribution in following way
\begin{eqnarray}
B_n(f;x)=\sum_{k=0}^{n}P_{n,k}(x)f\bigg(\frac{k}{n}\bigg), \hspace{2.0cm} n=1,2,3...,k=0,1,2,...
\end{eqnarray}
where $P_{n,k}(x)=(^n_k) x^k(1-x)^{n-k}$ and proved pointwise and uniform approximation in the space of continuous functions on $[0,1]$. These operators provide the powerful tool for numerical analysis, computer added geometric design(CAGD) and solutions of differential equations. But these operators are not suitable for discontinuous functions. Later on, Kantorovich\cite{Kantorovich} generalized the Bernstein operators  for integrable functions as
\begin{eqnarray}
  K_n(f;x)=(n+1)\sum\limits_{k=0}^{n}P_{n,k}(x)\int\limits_{\frac{k}{n+1}}^{\frac{k+1}{n+1}}f(t)dt, \hspace{2 cm} k=0,1,2,3,...,n=1,2,3,....
\end{eqnarray}
where $P_{n,k}(x)={n\choose k}x^k(1-x)^{n-k}$, $0\leq x<\leq 1$.  Szasz\cite{szasz} introduced linear positive operators in the sense of exponential growth on non-negative semi axes\\
\begin{eqnarray}
S_{n,k}(f;x)=\sum\limits_{k=0}^{\infty}s_{n,k}f\bigg(\frac{k}{n}\bigg), \hspace{2 cm} n=1,2,3,....
\end{eqnarray}
where $f(x)\in C[0,\infty)$ and $s_{n,k}=e^{-nx}\frac{(nx)^k}{k!}$. Several generalizations of these operators have been studied by different researchers (\cite{duman}-\cite{agarwal}). A generalization of operators (1) was given by Stancu\cite{D.D.} depending on the parameters $\alpha$ and $\beta$ such that $0\leq\alpha\leq\beta$ on $[0,1]$. Many operators preserve the constant and linear functions but these operators do not preserve $x^2$. King\cite{king} introduced a method in order to preserve $x^2$ for the Bernstein operators.

Recently, Varma and Tasdelen\cite{Varma} gave a generalization of well known Szasz-Mirakjan operators using Charlier polynomials\cite{ismail} having the generating function of the form
\begin{eqnarray}
 e^t\bigg(1-\frac{t}{a}\bigg)^u=\sum_{k=0}^{\infty} C^{(a)}_k(u)\frac{t^k}{k!}, \hspace{2cm} |t|<a
\end{eqnarray}
and the explicit representation
\begin{eqnarray*}
C_k^{(a)}(u)=\sum\limits_{r=0}^{k}{k \choose r}(-u)_r\Bigg(\frac{1}{a}\Bigg)^r,
\end{eqnarray*}
where $(\alpha)_k$ is the Pochhammer's symbol given by\\
\begin{eqnarray*}
(\alpha)_0=1, (\alpha)_k=\alpha(\alpha+1)...(\alpha+k-1),  k=1,2,....
\end{eqnarray*}
 We note that for $a>0$ and $u\leq 0$, Charlier polynomials are positive. Varma and Tadelen\cite{Varma} defined the Szasz-type and Kantorovich-Szasz-type operators operators as
\begin{eqnarray}
L_n(f;x,a)=e^{-1}\bigg(1-\frac{1}{a}\bigg)^{(a-1)nx}\sum_{k=0}^{\infty}\frac{C^{(u)}_k(-(a-1)nx)}{k!}f\bigg(\frac{k}{n}\bigg),\\
L^{*}_n(f;x,a)=ne^{-1}\bigg(1-\frac{1}{a}\bigg)^{(a-1)nx}\sum_{k=0}^{\infty}\frac{C^{(u)}_k(-(a-1)nx)}{k!}\int\limits_{\frac{k}{n}}^{\frac{k+1}{n}}f(s)ds,
\end{eqnarray}
 where $a>1$, $n=1,2,3...,k=0,1,2...$ and $x\geq 0$. \\
 Recently, Wafi and Rao\cite{rao} generalized the operators (5) as follows:
 \begin{eqnarray}
 T_{n,a}^{\alpha,\beta}(f;r_{n,a}(x;\alpha,\beta))= e^{-1}\bigg(1-\frac{1}{a}\bigg)^{(a-1)nr_{n,a}(x;\alpha,\beta)}\sum_{k=0}^{\infty}\frac{C^{(u)}_k(-(a-1)nr_{n,a}(x;\alpha,\beta))}{k!}f\bigg(\frac{k+\alpha}{n+\beta}\bigg)
 \end{eqnarray}
 for any function $f\in C[0,\infty)$, $x \geq 0$ and $0\leq\alpha\leq\beta$. Where
 \begin{eqnarray}
  r_{n,a}(x;\alpha,\beta)=\frac{-(3+2\alpha+\frac{1}{a-1})+\sqrt{(3+2\alpha+\frac{1}{a-1})^2+4((n+\beta)^2x^2-2-2\alpha-\alpha^2)}}{2n}.
 \end{eqnarray}
   In this paper, they discussed the rate of convergence and obtained better error estimation than (5). Motivated by the above development, we define a Kantorovich version of $T_{n,a}^{\alpha,\beta}$ by
\begin{eqnarray}
 K_{n,a}^{\alpha,\beta}(f;r^{*}_{n,a}(x;\alpha,\beta),a)=(n+\beta)e^{-1}\bigg(1-\frac{1}{a}\bigg)^{(a-1)nr^{*}_{n,a}(x;\alpha,\beta)}\sum_{k=0}^{\infty}\frac{C^{(u)}_k(-(a-1)nr^{*}_{n,a}(x;\alpha,\beta))}{k!}\int\limits_{\frac{k+\alpha}{(n+\beta)}}^{\frac{k+\alpha+1}{(n+\beta)}}f(s)ds,
 \end{eqnarray}
 where
 \begin{eqnarray}
  r^{*}_{n,a}(x;\alpha,\beta)=\frac{-(4+2\alpha+\frac{1}{a-1})+\sqrt{(4+2\alpha+\frac{1}{a-1})^2+4((n+\beta)^2x^2-\frac{10}{3}-3\alpha-\alpha^2)}}{2n}.
 \end{eqnarray}
 \vspace{5 cm}\\
We observe that\\ \\
 (i) if $\alpha=\beta=0$ and $r^{*}_{n,a}(x;\alpha,\beta)=x$, operators (9) reduce to operators (6),\\
 and\\
 (ii) for $\alpha=\beta=0$ and $r^{*}_{n,a}(x;\alpha,\beta)=x$ as $a\rightarrow \infty$ and taking $x-\frac{1}{n}$ instead of $x$, operators
 (9) reduce to the Classical Kantorovich-Szasz operators.\\

 In the present paper, we discuss the rate of convergence for continuous functions, first order derivative of the function   and weighted Korovkin type theorem. Further, we investigate some direct and local approximation results using Ditzian-Totik modulus of smoothness, second order modulus of continuity, Peetre's K-functional and Lipschitz space.
 \section{Basic Estimates}
 \textbf{Lemma 2.1} From Wafi and Rao\cite{rao}, we have
 \begin{eqnarray*}
\sum\limits_{k=0}^{\infty}\frac{C_k^{(u)}(-(a-1)nr_{n,a}(x;\alpha,\beta))}{k!}&=&e\Bigg(1-\frac{1}{a}\Bigg)^{-(a-1)nr_{n,a}(x;\alpha,\beta)},\\
\sum\limits_{k=0}^{\infty}k\frac{C_k^{(u)}(-(a-1)nr_{n,a}(x;\alpha,\beta))}{k!}&=&e\Bigg(1-\frac{1}{a}\Bigg)^{-(a-1)nr_{n,a}(x;\alpha,\beta)}(1+nr_{n,a}(x;\alpha,\beta)),\\
\nonumber\sum\limits_{k=0}^{\infty}k^2\frac{C_k^{(u)}(-(a-1)nr_{n,a}(x;\alpha,\beta))}{k!}&=&e\Bigg(1-\frac{1}{a}\Bigg)^{-(a-1)nr_{n,a}(x;\alpha,\beta)}\Bigg(2+\bigg(3+\frac{1}{a-1}\bigg)nr_{n,a}(x;\alpha,\beta)\\
&&+n^2r_n^2(x;\alpha,\beta)(x)\Bigg),\\
\sum\limits_{k=0}^{\infty}k^3\frac{C_k^{(u)}(-(a-1)nr_{n,a}(x;\alpha,\beta))}{k!}&=&e\Bigg(1-\frac{1}{a}\Bigg)^{-(a-1)nr_{n,a}(x;\alpha,\beta)}\Bigg(6+\bigg(11+\frac{6}{a-1}+\frac{2}{(a-1)^2}\bigg)nr_{n,a}(x;\alpha,\beta)\\
&&+\bigg(6+\frac{3}{a-1}\bigg)n^2r_n^2(x;\alpha,\beta)+n^3r_{n,a}^3(x;\alpha,\beta)\Bigg),\\   \sum\limits_{k=0}^{\infty}k^4\frac{C_k^{(u)}(-(a-1)nr_{n,a}(x;\alpha,\beta))}{k!}&=&e\Bigg(1-\frac{1}{a}\Bigg)^{-(a-1)nr_{n,a}(x;\alpha,\beta)}\Bigg(24+\bigg(50+\frac{35}{a-1}+\frac{20}{(a-1)^2}+\frac{6}{(a-1)^3}\bigg)\\
&&\times nr_{n,a}(x;\alpha,\beta)+\bigg(35+\frac{30}{a-1}+\frac{11}{a-1}^2\bigg)n^2r_{n,a}^2(x;\alpha,\beta)\\
&&+\bigg(10+\frac{6}{a-1}\bigg)n^3r_{n,a}^3(x;\alpha,\beta)+n^4r_{n,a}^4(x;\alpha,\beta)\Bigg).
 \end{eqnarray*}
\textbf{Lemma 2.2} Let $e_i=t^i,i=0,1,2$. Then for the operators $K_{n,a}^{\alpha,\beta}$, we have
\begin{eqnarray*}
&(i)\quad K_{n,a}^{\alpha,\beta}(1;x)=&1\\
&(ii)\quad K_{n,a}^{\alpha,\beta}(t;x)=&\frac{-(1+\frac{1}{a-1})+\sqrt{-(4+2\alpha+\frac{1}{a-1})^2+4((n+\beta)^2x^2-\frac{10}{3}-3\alpha-\alpha^2)}}{2(n+\beta)}.\\
&(iii)\quad K_{n,a}^{\alpha,\beta}(t^2;x)=&x^2.
\end{eqnarray*}
\textbf{Proof} Replacing $r_{n,a}(x;\alpha,\beta)$ by $r^{*}_{n,a}(x;\alpha,\beta)$ in the Lemma 2.1,  we have
\begin{eqnarray*}
(i)\quad K_{n,a}^{\alpha,\beta}(1;x,a)&=& (n+\beta)e^{-1}\bigg(1-\frac{1}{a}\bigg)^{(a-1)nr^{*}_{n,a}(x;\alpha,\beta)}\sum_{k=0}^{\infty} \frac{C^{(u)}_k(-(a-1)nr^{*}_{n,a}(x;\alpha,\beta))}{k!}\int\limits_{\frac{k+\alpha}{n+\beta}}^{\frac{k+\alpha+1}{n+\beta}}1.dt\quad\quad\\
&=&(n+\beta)e^{-1}\bigg(1-\frac{1}{a}\bigg)^{(a-1)nr^{*}_{n,a}(x;\alpha,\beta)}\sum_{k=0}^{\infty} \frac{C^{(u)}_k(-(a-1)nr^{*}_{n,a}(x;\alpha,\beta))}{k!}\bigg[\frac{k+\alpha+1}{n+\beta}-\frac{k+\alpha}{n+\beta}\bigg]\\
&&=(n+\beta)e^{-1}\bigg(1-\frac{1}{a}\bigg)^{(a-1)nr^{*}_{n,a}(x;\alpha,\beta)}e\bigg(1-\frac{1}{a}\bigg)^{-(a-1)nr^{*}_{n,a}(x;\alpha,\beta)}\times\frac{1}{n+\beta}\\
&&=1,\\
(ii)\quad K_{n,a}^{\alpha,\beta}(t;x,a)&=& (n+\beta) e^{-1}\bigg(1-\frac{1}{a}\bigg)^{(a-1)nr^{*}_{n,a}(x;\alpha,\beta)}\sum_{k=0}^{\infty}\frac{C^{(u)}_k(-(a-1)nr^{*}_{n,a}(x;\alpha,\beta))}{k!}\int\limits_{\frac{k+\alpha}{n+\beta}}^{\frac{k+\alpha+1}{n+\beta}}t.dt\\
&&=(n+\beta)e^{-1}\bigg(1-\frac{1}{a}\bigg)^{(a-1)nr^{*}_{n,a}(x;\alpha,\beta)}\sum_{k=0}^{\infty}\frac{C^{(u)}_k(-(a-1)nr^{*}_{n,a}(x;\alpha,\beta))}{k!}\frac{1}{2}\bigg[\bigg(\frac{k+\alpha+1}{n+\beta}\bigg)^2-\bigg(\frac{k+\alpha}{n+\beta}\bigg)^2\bigg]\\
&&=(n+\beta)e^{-1}\bigg(1-\frac{1}{a}\bigg)^{(a-1)nr^{*}_{n,a}(x;\alpha,\beta)}\sum_{k=0}^{\infty}\frac{C^{(u)}_k(-(a-1)nr^{*}_{n,a}(x;\alpha,\beta))}{k!}\frac{2(k+\alpha)+1}{2(n+\beta)^2}\\
&&=e^{-1}\bigg(1-\frac{1}{a}\bigg)^{(a-1)nr^{*}_{n,a}(x;\alpha,\beta)}\sum_{k=0}^{\infty}\frac{C^{(u)}_k(-(a-1)nr^{*}_{n,a}(x;\alpha,\beta))}{k!}\frac{(k+\alpha)}{(n+\beta)}+\frac{1}{2(n+\beta)}.\\
&&=\frac{1}{n+\beta}e^{-1}\bigg(1-\frac{1}{a}\bigg)^{(a-1)nr^{*}_{n,a}(x;\alpha,\beta)}\sum_{k=0}^{\infty}\frac{C^{(u)}_k(-(a-1)nr^{*}_{n,a}(x;\alpha,\beta))}{k!}k+\frac{2\alpha+1}{2(n+\beta)}\\
&&=\frac{1}{n+\beta}e^{-1}\bigg(1-\frac{1}{a}\bigg)^{(a-1)nr^{*}_{n,a}(x;\alpha,\beta)}e\Bigg(1-\frac{1}{a}\Bigg)^{-(a-1)nr^{*}_{n,a}(x;\alpha,\beta)}(1+nr^{*}_{n,a}(x;\alpha,\beta))+\frac{2\alpha+1}{2(n+\beta)}\\
&&=\frac{nr^{*}_{n,a}(x;\alpha,\beta)}{2(n+\beta)}+\frac{2\alpha+3}{2(n+\beta)}\\
&&=\frac{-(4+2\alpha+\frac{1}{a-1})+\sqrt{(4+2\alpha+\frac{1}{a-1})^2+4((n+\beta)^2x^2-\frac{10}{3}-3\alpha-\alpha^2)}}{2(n+\beta)}+\frac{2\alpha+3}{2(n+\beta)}\\
&&=\frac{-(1+\frac{1}{a-1})+\sqrt{(4+2\alpha+\frac{1}{a-1})^2+4((n+\beta)^2x^2-\frac{10}{3}-3\alpha-\alpha^2)}}{2(n+\beta)}
\end{eqnarray*}
Similarly, we can prove
\begin{eqnarray*}
 K_{n,a}^{\alpha,\beta}(t^2;x,a)&=& x^2.
\end{eqnarray*}
\textbf{Lemma 2.3} Let $\psi_x^i(t)=(t-x)^i,i=0,1,2$. Then
\begin{eqnarray*}
K_{n,a}^{\alpha,\beta}(\psi_x^0;x)&=&1,\\
K_{n,a}^{\alpha,\beta}(\psi_x^1;x)&=&-\frac{\big(1+\frac{1}{a-1}\big)}{2(n+\beta)}+\frac{\frac{8}{3}+4\alpha+\frac{4\alpha+8}{a-1}}{(n+\beta)\big(\sqrt{(4+2\alpha+\frac{1}{a-1})^2+4((n+\beta)^2x^2-\frac{10}{3}-3\alpha-\alpha^2)}+2(n+\beta)x\big)},\\
K_{n,a}^{\alpha,\beta}(\psi_x^2;x)&=&\frac{(1+\frac{1}{a-1})}{n+\beta}x-\frac{x}{(n+\beta)}\frac{\big(\frac{8}{3}+4\alpha+\frac{4\alpha+8}{a-1}\big)}{\sqrt{(4+2\alpha+\frac{1}{a-1})^2+4((n+\beta)^2x^2-\frac{10}{3}-3\alpha-\alpha^2)}+2(n+\beta)x}. \end{eqnarray*}
\textbf{Proof} In view of Lemma 2.2 and linearity property we can easily prove this Lemma.\\ \\
\section{Order of approximation for the function $f$ and derivative of $f$ }
Let $f\in C[0,\infty)$. Then modulus of continuity od $f$ defined as follows
\begin{eqnarray*}
\omega(f;\delta)=\sup_{|t-y|\leq \delta}|f(t)-f(y)|, \hspace{1 cm} t,y\in[0,\infty).
\end{eqnarray*}
For $f\in C[0,\infty)$ and $\delta>0$, one has
\begin{eqnarray}
 |f(t)-f(y)|\leq \bigg(1+\frac{(t-y)^2}{\delta^2}\bigg)\omega(f;\delta).
\end{eqnarray}
And \\ \\
$ E=\{f:[0,\infty)\rightarrow R, |f(x)|\leq Me^{Ax}, A\in R$ and  $M\in R^+\}.$\\ \\
\textbf{Theorem 3.1} Let $f\in C[0,\infty)\cap E$ and $x\geq0$. Then for operators $K_{n,a}^{\alpha,\beta} $, we have
\begin{eqnarray*}
|K_{n,a}^{\alpha,\beta}(f;x)-f(x)|\leq 2\omega\bigg(f;\delta_{n,a}^{\alpha,\beta}\bigg),
\end{eqnarray*}
where $\delta_{n,a}^{\alpha,\beta}=\sqrt{K_{n,a}^{\alpha,\beta}(\psi_x^2;x)}$.\\ \\
\textbf{Proof} From (11), we have
\begin{eqnarray*}
 |K_{n,a}^{\alpha,\beta}(f;x)-f(x)|&\leq& (n+\beta) e^{-1}\Bigg(1-\frac{1}{a}\Bigg)^{(a-1)nr^{*}_{n,a}(x;\alpha,\beta)}\sum_{k=0}^{\infty}\frac{C^{(u)}_k(-(a-1)nr^{*}_{n,a}(x;\alpha,\beta))}{k!}\int\limits_{\frac{k+\alpha}{n+\beta}}^{\frac{k+\alpha+1}{n+\beta}}|f(t)-f(x)|dt\\
 &\leq&\Bigg\{(n+\beta) e^{-1}\Bigg(1-\frac{1}{a}\Bigg)^{(a-1)nr^{*}_{n,a}(x;\alpha,\beta)}\sum_{k=0}^{\infty}\frac{C^{(u)}_k(-(a-1)nr^{*}_{n,a}(x;\alpha,\beta))}{k!}\int\limits_{\frac{k+\alpha}{n+\beta}}^{\frac{k+\alpha+1}{n+\beta}}\bigg(1+\frac{(t-x)^2}{(\delta_{n,a}^{\alpha,\beta})^2}\bigg)dt\Bigg\}\\
 &&\times\omega(f;\delta_{n,a}^{\alpha,\beta})\\
 &\leq& \Bigg\{1+\frac{K_{n,a}^{\alpha,\beta}(\psi_x^2;x)}{(\delta_{n,a}^{\alpha,\beta})^2}\Bigg\}\omega(f;\delta_{n,a}^{\alpha,\beta})\\
&=& 2 \omega(f;\delta_{n,a}^{\alpha,\beta}).
\end{eqnarray*}
where $\delta_{n,a}^{\alpha,\beta}=\sqrt{K_{n,a}^{\alpha,\beta}(\psi_x^2;x)}$.\\ \\
\textbf{Remark} For the Kantorovich-Szasz type operators $L_n^*$ given by (6), we have, for every $f\in C[0,\infty)\cap E$,
\begin{eqnarray}
|L_n^*(f;x,a)-f(x)|\leq 2\omega(f;\delta),
\end{eqnarray}
where $\delta=\sqrt{\frac{x}{n}\bigg(1+\frac{1}{a-1}\bigg)+\frac{10}{3n^2}}.$
Here we show that our operators $K_{n,a}^{\alpha,\beta}$ has the better approximation than the operators $L_n^*$.\\
Since
\begin{eqnarray*}
\frac{x}{n+\beta}\bigg(1+\frac{1}{a-1}\bigg)<\frac{x}{n}\bigg(1+\frac{1}{a-1}\bigg)+\frac{10}{3n^2}.
\end{eqnarray*}
Then $\delta_{n,a}^{\alpha,\beta}<\delta$.\\ \\
 \textbf{Theorem 3.2} If $f'(x)$ has continuous derivative over $[0,\infty)$ and $\omega_1(f;\delta_{n,\beta})$ is the modulus of continuity of $f'(x)$, then, for $0\leq\alpha\leq\beta, a>1$ and $ x\in [0,b], b<\infty$, we have
\begin{eqnarray*}
 | K_{n,a}^{\alpha,\beta}(f;x)-f(x) | \leq\omega_1\big((n+\beta)^{-1}\big)\sqrt{K_{n,a}^{\alpha,\beta}(\psi_x^2(t);x)}\bigg\{1+\sqrt{(n+\beta)}\sqrt{K_{n,a}^{\alpha,\beta}(\psi_x^2(t);x)} \bigg\}.
\end{eqnarray*}
\textbf{Proof} It is known that
\begin{eqnarray}
\nonumber f(x_1)-f(x_2)&=&(x_1-x_2)f'(\xi),\\
 &=&(x_1-x_2)f'(x_1)+(x_1-x_2)[f'(\xi)-f'(x_1)],
\end{eqnarray}
for $x_1,x_2\in [0,b]$ and $x_1<\xi<x_2$. Also, we have
\begin{eqnarray}
 |(x_1-x_2)[f'(\xi)-f'(x_1)]| \leq| x_1-x_2 |(\lambda+1)\omega_1(\delta), \hspace{2cm} \lambda=\lambda(x_1,x_2;\delta).
\end{eqnarray}
Next, we find
\begin{eqnarray}
| K_{n,a}^{\alpha,\beta}(f;x)-f(x) |= \bigg| (n+\beta) e^{-1}\Bigg(1-\frac{1}{a}\Bigg)^{(a-1)nr^{*}_{n,a}(x;\alpha,\beta)}\sum_{k=0}^{\infty}\frac{C^{(u)}_k(-(a-1)nr^{*}_{n,a}(x;\alpha,\beta))}{k!}\int\limits_{\frac{k+\alpha}{n+\beta}}^{\frac{k+\alpha+1}{n+\beta}}f(t)-f(x)dt\bigg|.
\end{eqnarray}
Using (13) and (14), we get
\begin{eqnarray*}
| K_{n,a}^{\alpha,\beta}(f;x)-f(x)| &\leq & \Bigg|(n+\beta) e^{-1}\Bigg(1-\frac{1}{a}\Bigg)^{(a-1)nr^{*}_{n,a}(x;\alpha,\beta)}\sum_{k=0}^{\infty}\frac{C^{(u)}_k(-(a-1)nr^{*}_{n,a}(x;\alpha,\beta))}{k!}\int\limits_{\frac{k+\alpha}{n+\beta}}^{\frac{k+\alpha+1}{n+\beta}}(x-t)f'(x)dt\Bigg|\\
&&+ \omega_1(\delta_n^\beta)(\lambda+1)(n+\beta) e^{-1}\Bigg(1-\frac{1}{a}\Bigg)^{(a-1)nr^{*}_{n,a}(x;\alpha,\beta)}\sum_{k=0}^{\infty}\frac{C^{(u)}_k(-(a-1)nr^{*}_{n,a}(x;\alpha,\beta))}{k!}\int\limits_{\frac{k+\alpha}{n+\beta}}^{\frac{k+\alpha+1}{n+\beta}}|t-x|dt\\ &\leq & \omega_1(\delta_n^\beta)\Bigg\{e^{-1}\Bigg(1-\frac{1}{a}\Bigg)^{(a-1)nr^{*}_{n,a}(x;\alpha,\beta)}\sum_{k=0}^{\infty}\frac{C^{(u)}_k(-(a-1)nr^{*}_{n,a}(x;\alpha,\beta))}{k!}\int\limits_{\frac{k+\alpha}{n+\beta}}^{\frac{k+\alpha+1}{n+\beta}}|t-x|dt\\
&+&e^{-1}\Bigg(1-\frac{1}{a}\Bigg)^{(a-1)nr^{*}_{n,a}(x;\alpha,\beta)}\sum_{\lambda\geq1}C^{(u)}_k(-(a-1)nr^{*}_{n,a}(x;\alpha,\beta))\int\limits_{\frac{k+\alpha}{n+\beta}}^{\frac{k+\alpha+1}{n+\beta}}|t-x|\lambda(x,t;\delta)dt\Bigg\}\\
&\leq & \omega_1(\delta_n^\beta)\Bigg\{e^{-1}\Bigg(1-\frac{1}{a}\Bigg)^{(a-1)nr^{*}_{n,a}(x;\alpha,\beta)}\sum_{k=0}^{\infty}\frac{C^{(u)}_k(-(a-1)nr^{*}_{n,a}(x;\alpha,\beta))}{k!}\int\limits_{\frac{k+\alpha}{n+\beta}}^{\frac{k+\alpha+1}{n+\beta}}|t-x|dt\\
&+&\frac{1}{\delta_n^\beta}e^{-1}\Bigg(1-\frac{1}{a}\Bigg)^{(a-1)nr^{*}_{n,a}(x;\alpha,\beta)}\sum_{k=0}^{\infty}\frac{C^{(u)}_k(-(a-1)nr^{*}_{n,a}(x;\alpha,\beta))}{k!}\int\limits_{\frac{k+\alpha}{n+\beta}}^{\frac{k+\alpha+1}{n+\beta}}(t-x)^2dt\Bigg\}\\
&\leq&\omega_1(\delta_n^\beta)\Bigg(\sqrt{K_{n,a}^{\alpha,\beta}(\psi_x^2;x)}+\frac{K_{n,a}^{\alpha,\beta}(\psi_x^2;x)}{\delta_n^\beta}\Bigg)\\
&=&\omega_1(\delta_n^\beta)\sqrt{K_{n,a}^{\alpha,\beta}(\psi_x^2;x)}\Bigg\{1+\frac{\sqrt{K_{n,a}^{\alpha,\beta}(\psi_x^2;x)}}{\delta_n^\beta}\Bigg\}
\end{eqnarray*}
Taking $\delta_n^\beta=(n+\beta)^{-1}$, we get
\begin{eqnarray*}
 | K_{n,a}^{\alpha,\beta}(f;x)-f(x) | \leq\omega_1((n+\beta)^{-1})\sqrt{K_{n,a}^{\alpha,\beta}(\psi_x^2(t);x)}\bigg\{1+\sqrt{(n+\beta)}\sqrt{K_{n,a}^{\alpha,\beta}(\psi_x^2(t);x)} \bigg\}.
\end{eqnarray*}

We shall now discuss the Korovkin-type theorem in polynomial weighted space of continuous and unbounded functions defined on $[0,\infty)$.   Here we recall some symbols and notions from\cite{Gad2}. Let $\rho(x)=1+x^2$, $-\infty<x<\infty$ and $B_\rho[0,\infty)=\{f(x):|f(x)|\leq M_f \rho(x),\rho(x)$ is weight function, $M_f$  is a constant depending on $f$ and $x\in[0,\infty) \}$. $C_\rho[0,\infty)$ is the space of continuous function in $B_\rho[0,\infty)$ with the norm $\|f(x)\|_\rho=\sup\limits_{x\in[0,\infty)}\frac{|f(x)|}{\rho(x)}$ and \\ $C_\rho^{k}=\{f\in C_\rho: \lim\limits_{|x|\rightarrow\infty}\frac{f(x)}{\rho(x)}=k,$ where $k$ is a constant depending on $f\}$.\\ \\
 \textbf{Theorem 3.3} Let $f\in C_\rho^k$. Then, the operators $ K_{n,a}^{\alpha,\beta}(f;x)$ converges uniformly to $f(x)$ on $[0,\infty)$.\\ \\
 \textbf{Proof} By the Korovkin Theorem, it is sufficient to verify that \\
 \begin{eqnarray*}
 \lim\limits_{n\rightarrow\infty}||K_{n,a}^{\alpha,\beta}(t^i;x)-x^i||_\rho=0,\hspace{1cm}for\hspace{0.3cm} i=0,1,2.
 \end{eqnarray*}
 It is obvious that  $\lim\limits_{n\rightarrow\infty}\|K_{n,a}^{\alpha,\beta}(1;x)-1\|_\rho=0$ and $\lim\limits_{n\rightarrow\infty}\|K_{n,a}^{\alpha,\beta}(x^2;x)-x^2\|_\rho=0$. Now, from the Lemma 2.1
  \begin{eqnarray*}
 \sup\limits_{x\in[0,\infty)}\frac{| K_{n,a}^{\alpha,\beta}(t;x)-x|}{1+x^2}&=&\sup\limits_{x\in[0,\infty)}\frac{\Bigg|\frac{-(1+\frac{1}{a-1})+\sqrt{\big(4+2\alpha+\frac{1}{a-1}\big)^2+4\big((n+\beta)^2x^2-\frac{10}{3}-3\alpha-\alpha^2\big)}}{2(n+\beta)}-x\Bigg|}{1+x^2}\\
 &=&\sup\limits_{x\in[0,\infty)}\frac{\Bigg|-\frac{\big(1+\frac{1}{a-1}\big)}{2(n+\beta)}+\frac{\frac{8}{3}+4\alpha+\frac{4\alpha+8}{a-1}}{2(n+\beta)\sqrt{\big(4+2\alpha+\frac{1}{a-1}\big)^2+4((n+\beta)^2x^2-\frac{10}{3}-3\alpha-\alpha^2)}+2(n+\beta)x}\Bigg|}{1+x^2}\\
 &\leq&\sup\limits_{x\in[0,\infty)}\Bigg|\frac{\frac{4}{3}+2\alpha+\frac{2\alpha+4}{a-1}}{(n+\beta)\sqrt{\big(4+2\alpha+\frac{1}{a-1}\big)^2+4\big((n+\beta)^2x^2-\frac{10}{3}-3\alpha-\alpha^2\big)}+2(n+\beta)x}\frac{1}{(1+x^2)}\Bigg|\\
 &&+\sup\limits_{x\in[0,\infty)}\Bigg|\frac{\big(1+\frac{1}{a-1}\big)}{2(n+\beta)}\frac{1}{(1+x^2)}\Bigg|
 \end{eqnarray*}
 which shows that $|| K_{n,a}^{\alpha,\beta}(t;x)-x||_\rho\rightarrow 0$. as $n\rightarrow\infty$,\\
 Hence, we proved the theorem.\\ \\
\section{Direct Estimate}
Ditzian-Totik Modulus of smoothness\cite{dit1} is defined as:
\begin{eqnarray*}
\omega^2_{\varphi^\lambda}(f;\delta)&=& \sup\limits_{0<h\leq \delta}\parallel \Delta^2_{h\varphi(x)}f(x) \parallel ,\\
&=& \sup\limits_{0<h\leq \delta} \hspace{0.5cm} \sup\limits_{x\pm h \varphi^\lambda \in [0,\infty)} |f(x-h\varphi^\lambda (x))-2f(x)+f(x+h\varphi^\lambda (x))|,
\end{eqnarray*}
where $\varphi^2(x)=x$.
And, Peetre's K-functional\cite{dit1} is given by
\begin{eqnarray}
K_{\varphi^\lambda}(f,\delta^2)=\inf\limits_g \bigg(\|f-g\|_{C[0,\infty)} +\delta^2\|\varphi^2\lambda g''\|_{C[0,\infty)}\bigg), \hspace{2cm} g,g'\in AC_{loc}.
\end{eqnarray}
The K-functional is equivalent to the modulus of smoothness, i.e.,
\begin{eqnarray}
C^{-1}K_{\varphi^\lambda}(f,\delta^2)\leq \omega^2_{\varphi^\lambda}(f,\delta)\leq C K_{\varphi^\lambda}(f,\delta^2).
\end{eqnarray}
First result based on Ditziaz-Totik modulus of smoothness was given by Ditzian\cite{dit} for the Bernstein polynomials as:
\begin{eqnarray*}
| B_n(f;x)-f(x)| \leq C \omega^2_{\varphi^\lambda}(f,n^{-\frac{1}{2}} \varphi(x)^{1-\lambda}).
\end{eqnarray*}
\vspace{4 cm}

Now, we prove the similar result for the operator $K_{n,a}^{\alpha,\beta}$. \newline \\
\textbf{Theorem 4.1} For $f\in L_p(0,\infty), 0\leq p<\infty$, and $a>1$, we have
\begin{eqnarray*}
 | K_{n,a}^{\alpha,\beta}(f;x)-f(x)| \leq C \omega^2_{\varphi^\lambda}\big(f,(n+\beta)^{-\frac{1}{2}}  \varphi(x)^{1-\lambda}\big)\hspace{0.2cm} for\hspace{0.2cm} large \hspace{0.2cm}n\\
 \end{eqnarray*}
 where $0\leq\lambda\leq 1$, $\varphi^2(x)=x$.\\ \\
 \textbf{Proof} Using (14),(15), we have
 \begin{eqnarray}
 \parallel\ f-g \parallel_{L_p(0,\infty)} \leq A\omega^2_{\varphi^\lambda}\big(f,(n+\beta)^{-\frac{1}{2}}\varphi(x)^{1-\lambda}\big),\\
 (n+\beta)^{-1}\varphi(x)^{2-2\lambda}\| \varphi^{2\lambda}g''\|_{L_p(0,\infty)}\leq B\omega^2_{\varphi^\lambda}\big(f,(n+\beta)^{-\frac{1}{2}}\varphi(x)^{1-\lambda}\big).
\end{eqnarray}
Next, we can choose $g_n\equiv g_{n,x,\lambda}$ for fixed $x$ and $\lambda+1$ such that
\begin{eqnarray*}
 | K_{n,a}^{\alpha,\beta}(f;x)-f(x)| & \leq & | K_{n,a}^{\alpha,\beta}(f-g_n;x)-(f-g_n)(x) | +| K_{n,a}^{\alpha,\beta}(g_n;x)-g_n(x)|,\\ \\
 &\leq & 2\parallel f-g_n \parallel_{L_p(0,\infty)}+| K_{n,a}^{\alpha,\beta}(g_n;x)-g_n(x)|.
\end{eqnarray*}
From (16), we get
\begin{eqnarray}
 | K_{n,a}^{\alpha,\beta}(f;x)-f(x)| & \leq 2A\omega^2_{\varphi^\lambda}\big(f,(n+\beta)^{-\frac{1}{2}}\varphi(x)^{1-\lambda}\big) +| K_{n,a}^{\alpha,\beta}(g_n;x)-g_n(x)|.
\end{eqnarray}
Now, the last term can be calculated by using Taylor's formula
 \begin{eqnarray*}
 | K_{n,a}^{\alpha,\beta}(g_n(t)-g_n(x);x)| &\leq& | g'_n(x)K_{n,a}^{\alpha,\beta}((t-x);x)|+\Big| K_{n,a}^{\alpha,\beta}\bigg(\int\limits_{t}^{x}(x-u)g''_n(u)du;x\bigg)\Big|\\
 &\leq& K_{n,a}^{\alpha,\beta}\bigg( \frac{| x-\frac{k}{n}|}{\varphi^{2\lambda}(x)}\int\limits_{\frac{k}{n}}^x \varphi^{2\lambda}(u)|g''_n(u)du|;x \bigg)\\
 &\leq & \parallel \varphi^{2\lambda}g''_n\parallel_{L_p(0,\infty)}\frac{1}{\varphi^{2\lambda}(x)}K_{n,a}^{\alpha,\beta}((t-x)^2;x)\\
 &\leq &  \parallel \varphi^{2\lambda}g''_n\parallel_{L_p(0,\infty)}\frac{1}{\varphi^{2\lambda}(x)}\frac{x}{n+\beta}\frac{(n+\beta)K_{n,a}^{\alpha,\beta}((t-x)^2;x)}{x}\\
 &\leq &  \parallel \varphi^{2\lambda}g''_n\parallel_{L_p(0,\infty)}\frac{x(n+\beta)^{-1}}{\varphi^{2\lambda}(x)}\frac{(n+\beta)K_{n,a}^{\alpha,\beta}((t-x)^2;x)}{x}.\\
 \end{eqnarray*}
For the large value of $n$, we get
\begin{eqnarray*}
\frac{(n+\beta)K_{n,a}^{\alpha,\beta}((t-x)^2;x)}{x}\leq \bigg(1+\frac{1}{a-1}\bigg).
\end{eqnarray*}
 Therefore
 \begin{eqnarray}
| K_{n,a}^{\alpha,\beta}(g_n(t)-g_n(x);x)| &\leq& \bigg(1+\frac{1}{a-1}\bigg)B\omega^2_{\varphi^\lambda}\big(f,(n+\beta)^{-\frac{1}{2}}\varphi(x)^{1-\lambda}.
 \end{eqnarray}
Using (18) and (19), we get
\begin{eqnarray*}
 | K_{n,a}^{\alpha,\beta}(f(t)-f(x);x)| &\leq& M\omega^2_{\lambda}\bigg(f,(n+\beta)^{\frac{-1}{2}}\varphi(x)^{1-\lambda}\bigg)
\end{eqnarray*}
where $M=max\Bigg(2A,\bigg(1+\frac{1}{a-1}\bigg)B\Bigg)$.\\ \\

Let $C_B[0,\infty)$ denote the space of real valued continuous and bounded functions $f$ on $[0,\infty)$ endowed with the norm
\begin{eqnarray*}
\|f\|=\sup\limits_{0\leq x<\infty}|f(x)|.
\end{eqnarray*}
Then, for any $\delta>0$, Peeter's K-functional is defined as
\begin{eqnarray*}
K_2(f,\delta)=inf\{\|f-g\|+\delta\|g''\|: g\in C_B^2[0,\infty)\},
\end{eqnarray*}
where $C_B^2[0,\infty)=\{g\in C_B[0,\infty):g',g''\in C_B[0,\infty)\}$. By Devore and Lorentz[\cite{devor}, p.177, Theorem 2.4], there exits an absolute constant $C>0$ such that
\begin{eqnarray*}
K_2(f;\delta)\leq C\omega_2(f;\sqrt{\delta}),
\end{eqnarray*}
where $\omega_2(f;\delta)$ is the second order modulus of continuity is defined as
\begin{eqnarray*}
\omega_2(f,\sqrt{\delta})=\sup\limits_{0<h\leq\sqrt{\delta}} \sup\limits_{x\in[0,\infty)} |f(x+2h)-2f(x+h)+f(x)|.
\end{eqnarray*}
\textbf{Theorem 4.2} Let $f\in C_B^2[0,\infty)$. Then for all $x\in[0,\infty)$ there exist a constant $K>0$ such that
\begin{eqnarray*}
\mid K_{n,a}^{\alpha,\beta}(f;x)-f(x)\mid\leq K\omega_2(f;\sqrt{\Pi_{n,a}^{\alpha,\beta}(x)})+\omega(f;\Lambda_{n,a}^{\alpha,\beta})
\end{eqnarray*}
where $\Lambda_{n,a}^{\alpha,\beta}(x)=K_{n,a}^{\alpha,\beta}(\psi_x;x)$ and  $\Pi_{n,a}^{\alpha,\beta}(x)=K_{n,a}^{\alpha,\beta}(\psi_x^2;x)+\big(K_{n,a}^{\alpha,\beta}(\psi_x;x)\big)^2.$\\ \\
\textbf{Proof} First, we define the auxiliary operators
\begin{eqnarray}
\widehat{K}_{n,a}^{\alpha,\beta}(f;x)=K_{n,a}^{\alpha,\beta}(f;x)+f(x)-f\big(\Lambda_{n,a}^{\alpha,\beta}(x)\big)
\end{eqnarray}
We find that
\begin{eqnarray*}
\widehat{K}_{n,a}^{\alpha,\beta}(1;x)=1,
\end{eqnarray*}
\begin{eqnarray*}
\widehat{K}_{n,a}^{\alpha,\beta}(\psi_x(t);x)=0
\end{eqnarray*}
\begin{eqnarray}
|\widehat{K}_{n,a}^{\alpha,\beta}(f;x)|\leq 3\|f\|.
\end{eqnarray}
Let $g\in C_B^2[0,\infty)$. By the Taylor's theorem\\
\begin{eqnarray}
g(t)=g(x)+(t-x)g'(x)+\int\limits_x^t (t-v)g''(v)dv,
\end{eqnarray}
Now
\begin{eqnarray*}
\widehat{K}_{n,a}^{\alpha,\beta}(g;x)-g(x)&=&g'(x)\widehat{K}_{n,a}^{\alpha,\beta}(t-x;x)+\widehat{K}_{n,a}^{\alpha,\beta}\Big( \int\limits_x^t (t-v)g''(v)dv;x\Big)\\
&=&\widehat{K}_{n,a}^{\alpha,\beta}\Big( \int\limits_x^t (t-v)g''(v)dv;x\Big)\\
&=&K_{n,a}^{\alpha,\beta}\Big( \int\limits_x^t (t-v)g''(v)dv;x\Big)-\int\limits_x^{\Lambda_{n,a}^{\alpha,\beta}}(\Lambda_{n,a}^{\alpha,\beta}-v)g''(v)dv.
\end{eqnarray*}
Therefore
\begin{eqnarray}
\nonumber | \widehat{K}_n^{\alpha,\beta}(g;x)-g(x)|&\leq&\Bigg|K_{n,a}^{\alpha,\beta}\Big( \int\limits_x^t (t-v)g''(v)dv;x\Big)\Bigg|+\Bigg|\int\limits_x^{\Lambda_{n,a}^{\alpha,\beta}} (\Lambda_{n,a}^{\alpha,\beta}-v)g''(v)dv\Bigg|.
\end{eqnarray}
Since
 \begin{eqnarray}
\Bigg| \int\limits_x^t (t-v)g''(v)dv\Bigg|\leq(t-x)^2\parallel g''\parallel
 \end{eqnarray}
 and
 \begin{eqnarray}
 \Bigg|\int\limits_x^{\Lambda_{n,a}^{\alpha,\beta}} (\Lambda_{n,a}^{\alpha,\beta}-v)g''(v)dv\Bigg|\leq(\Lambda_{n,a}^{\alpha,\beta})^2\parallel g''\parallel
 \end{eqnarray}
 Then from(24),(25)and (26) implies that
 \begin{eqnarray}
\nonumber | \widehat{K}_{n,a}^{\alpha,\beta}(g;x)-g(x)|&\leq& \Bigg\{K_{n,a}^{\alpha,\beta}((t-x)^2;x)+(\Lambda_{n,a}^{\alpha,\beta})^2\Bigg\}\|g''\|\\
 &=&\Pi_{n,a}^{\alpha,\beta}(x)\|g''\|
 \end{eqnarray}
 Next, we have
 \begin{eqnarray*}
 |K_{n,a}^{\alpha,\beta}(f;x)-f(x)|\leq |\widehat{K}_{n,a}^{\alpha,\beta}(f-g;x)|+|(f-g)(x)|+|\widehat{K}_{n,a}^{\alpha,\beta}(g;x)-g(x)|+\big|f(\Lambda_{n,a}^{\alpha,\beta})-f(x)\big|
 \end{eqnarray*}
 Using(27), we have
 \begin{eqnarray*}
 |K_{n,a}^{\alpha,\beta}(f;x)-f(x)|&\leq& 4\|f-g\| +\widehat{K}_{n,a}^{\alpha,\beta}(g;x)-g(x)|+\big|f(\Lambda_{n,a}^{\alpha,\beta})-f(x)\big|\\
 &\leq& 4\|f-g\|+\Pi_{n,a}^{\alpha,\beta}(x)\|g''\|+\omega\Big(f;\Lambda_{n,a}^{\alpha,\beta}\Big).
 \end{eqnarray*}
 By the definition of Peetre's K-functional
 \begin{eqnarray*}
 |K_{n,a}^{\alpha,\beta}(f;x)-f(x)|\leq C\omega_2\big(f;\sqrt{\gamma_n^{\alpha,\beta}(x)}\big)+\omega(f;\Lambda_{n,a}^{\alpha,\beta}).
 \end{eqnarray*}

  Here, we discuss a local result in Lipschitz class
 \begin{eqnarray*}
 Lip_M^*=\{f\in C[0,\infty):|f(t)-f(x)|\leq M\frac{|t-x|^\alpha}{(t+x)^{\frac{\alpha}{2}}}:x,t\in(0,\infty)\}
 \end{eqnarray*}
 where $M$ is a constant and $0<\alpha\leq 1$ to prove the following theorem:\\ \\ \\
 \textbf{Theorem 4.3} Let $f\in Lip_M^*(\alpha)$ and $x\in (0,\infty)$. Then, we have
 \begin{eqnarray*}
  |K_{n,a}^{\alpha,\beta}(f;x)-f(x)|\leq M\Bigg[\frac{\Theta_n^{\alpha,\beta}(x)}{x}\Bigg],
 \end{eqnarray*}
 where $\Theta_{n,a}^{\alpha,\beta}(x)=K_{n,a}^{\alpha,\beta}((t-x)^2;x)$. \\ \\
 \textbf{Proof} Let $\alpha=1$ and $x\in(0,\infty)$. Then, for $f\in Lip_M^*(1)$, we have
 \begin{eqnarray*}
  |K_{n,a}^{\alpha,\beta}(f;x)-f(x)|&\leq & (n+\beta)e^{-1}\Bigg(1-\frac{1}{a}\Bigg)^{(a-1)nr^{*}_{n,a}(x;\alpha,\beta)}\sum_{k=0}^{\infty}\frac{C^{(u)}_k(-(a-1)nr^{*}_{n,a}(x;\alpha,\beta))}{k!}\int_{\frac{k+\alpha}{n+\beta}}^{\frac{k+\alpha+1}{n+\beta}}|f(t)-f(x)|dt\\
  &\leq & M(n+\beta)e^{-1}\Bigg(1-\frac{1}{a}\Bigg)^{(a-1)nr^{*}_{n,a}(x;\alpha,\beta)}\sum_{k=0}^{\infty}\frac{C^{(u)}_k(-(a-1)nr^{*}_{n,a}(x;\alpha,\beta))}{k!}\int_{\frac{k+\alpha}{n+\beta}}^{\frac{k+\alpha+1}{n+\beta}}\frac{|t-x|}{\sqrt{t+x}}dt.\\
  &\leq & \frac{M}{\sqrt{x}}(n+\beta)e^{-1}\Bigg(1-\frac{1}{a}\Bigg)^{(a-1)nr^{*}_{n,a}(x;\alpha,\beta)}\sum_{k=0}^{\infty}\frac{C^{(u)}_k(-(a-1)nr^{*}_{n,a}(x;\alpha,\beta))}{k!}\int_{\frac{k+\alpha}{n+\beta}}^{\frac{k+\alpha+1}{n+\beta}}|t-x|dt\\
  &\leq&\frac{M}{\sqrt{x}}K_n^{\alpha,\beta}(|t-x|;x)\\
  &\leq& M \frac{\sqrt{K_n^{\alpha,\beta}((t-x)^2;x)}}{\sqrt{x}}\\
  &=&M \Bigg(\frac{\Theta_{n,a}^{\alpha,\beta}(x)}{x}\Bigg)^{\frac{1}{2}}
 \end{eqnarray*}
 Thus, the assertion hold for $\alpha=1$. Now, we will prove for $\alpha\in (0,1)$. From the Holder inequality with $p=\frac{1}{\alpha}, \\ q=\frac{1}{1-\alpha}$, we have
 \begin{eqnarray*}
  |K_{n,a}^{\alpha,\beta}(f;x)-f(x)|&=& \Bigg(e^{-1}\Bigg(1-\frac{1}{a}\Bigg)^{(a-1)nr^{*}_{n,a}(x;\alpha,\beta)}\sum_{k=0}^{\infty}\frac{C^{(u)}_k(-(a-1)nr^{*}_{n,a}(x;\alpha,\beta))}{k!}\Bigg((n+\beta)\int_{\frac{k+\alpha}{n+\beta}}^{\frac{k+\alpha+1}{n+\beta}}|f(t)-f(x)|dt\Bigg)^{\frac{1}{\alpha}}\Bigg)^\alpha\\
  &&\times\big(e^{-1}\Bigg(1-\frac{1}{a}\Bigg)^{(a-1)nr^{*}_{n,a}(x;\alpha,\beta)}\sum_{k=0}^{\infty}\frac{C^{(u)}_k(-(a-1)nr^{*}_{n,a}(x;\alpha,\beta))}{k!}\big)^{1-\alpha}\\
  &\leq& \Bigg(e^{-1}\Bigg(1-\frac{1}{a}\Bigg)^{(a-1)nr^{*}_{n,a}(x;\alpha,\beta)}\sum_{k=0}^{\infty}\frac{C^{(u)}_k(-(a-1)nr^{*}_{n,a}(x;\alpha,\beta))}{k!}\Bigg((n+\beta)\int_{\frac{k+\alpha}{n+\beta}}^{\frac{k+\alpha+1}{n+\beta}}|f(t)-f(x)|dt\Bigg)^{\frac{1}{\alpha}}\Bigg)^\alpha
 \end{eqnarray*}
 Since $f\in Lip_M^*$, we obtain
 \begin{eqnarray*}
  |K_{n,a}^{\alpha,\beta}(f;x)-f(x)|&\leq& M\Bigg((n+\beta)e^{-1}\Bigg(1-\frac{1}{a}\Bigg)^{(a-1)nr^{*}_{n,a}(x;\alpha,\beta)}\sum_{k=0}^{\infty}\frac{C^{(u)}_k(-(a-1)nr^{*}_{n,a}(x;\alpha,\beta))}{k!}\int_{\frac{k+\alpha}{n+\beta}}^{\frac{k+\alpha+1}{n+\beta}}\frac{|t-x|}{\sqrt{t+x}}dt\Bigg)^\alpha\\ &\leq&\frac{M}{x^{\frac{\alpha}{2}}}\Bigg((n+\beta)e^{-1}\Bigg(1-\frac{1}{a}\Bigg)^{(a-1)nr^{*}_{n,a}(x;\alpha,\beta)}\sum_{k=0}^{\infty}\frac{C^{(u)}_k(-(a-1)nr^{*}_{n,a}(x;\alpha,\beta))}{k!}\int_{\frac{k+\alpha}{n+\beta}}^{\frac{k+\alpha+1}{n+\beta}}|t-x|dt\Bigg)^\alpha\\
  &=&\frac{M}{x^{\frac{\alpha}{2}}}\big(K_{n,a}^{\alpha,\beta}(|t-x|;x)\big)^{\alpha}\\
  &\leq& M \Bigg(\frac{\Theta_{n,a}^{\alpha,\beta}(x)}{x}\Bigg)^{\frac{\alpha}{2}}.
  \end{eqnarray*}

\end{document}